# A model of information propagation in transportation networks

Omar Mansour [1]*, Tomer Toledo[1], Shadi Haj-Yahia[1], Wafa Elias[2]

[1] Faculty of Civil and Environmental Engineering, Technion – Israel Institute of Technology, Haifa 32000, Israel

[2] Department of Civil Engineering, Shamoon College of Engineering, 84 Jabotinsky St., 84100 Ashdod, Israel

## Abstract

This paper introduces a new macroscopic perspective for simulating transportation networks. The idea is to look at the network as connected nodes. Each node sends an "information package" to its neighbors. Basically, the information package contains a state change that a specific node experienced, and it might affect the traffic network state in the future. Different types of information can be counted in transportation network. Each information type has different characteristics. It propagates through the network and interacts with other IPs and nodes. As a result, the model enables implementing and analyzing dynamic and inconvenient control strategies.

This paper focus on flow dynamics and demand routing information under complex environment. The flow dynamics flows the LWR theory along the links. The demand routing follows a node equilibrium model. The node model takes into account the users' preferences and choices under dynamic control agents such as dynamic tolling and sends the needed information for other facilities to operate.

A DNL model is developed based on the IPM. Several case studies demonstrate the use IPM and its potential to provide a reliable and results for real world complicated transportation applications.

## 1. Introduction

Dynamic network loading (DNL) models express the propagation of trip demands in the network using traffic flow principles. They are central to dynamic traffic assignment (DTA) models, which are widely used to evaluate urban transportation systems. With technology advancements, DNLs are also used to estimate and predict traffic conditions in real-time to make appropriate control actions (Toledo et al., 2017). In general, DNL models may be classified to three main scale groups: microscopic, mesoscopic and macroscopic. These scales represent the level of detail in modeling the flow entities and their movement in the network. Macroscopic models, which treat traffic flow as a continuous stream, are the least detailed but provide greater computational efficiency.

Macroscopic DNL models differ in the way they describe traffic dynamics within links. The LWR kinematic wave theory (Lighthill et al., 1955; Richards, 1956) describes traffic dynamics using a system of partial differential equations (PDE) that represent flow conservation and a fundamental diagram (FD):

$$\frac{\partial k(x,t)}{\partial t} + \frac{\partial f(x,t)}{\partial x} = 0 \tag{1}$$

$$f(k) = v(k)k \tag{2}$$

Where, $k$ is density, $f$ is flow, $x$ is a space location and $t$ a time point. $v$ is the space-mean speed.

The system above constitutes a first order model, which allows instantaneous speed changes. Second order models add equations to describe smooth changes in speed when traffic densities change (Messmer et al., 1990). The cell transmission model (CTM) (Daganzo, 1994) and its variants (e.g., Gomes et al., 2006; Srivastava et al., 2015) solve this formulation. Assuming a triangular FD, the CTM applies a finite differences solution scheme that discretizes both time and space. This discretization introduces approximation errors, which propagate in the network and affect the solution's accuracy (Raadsen et al., 2016).

The Link Transmission Model (LTM) family of models (Himpe et al., 2016; Newell, 1993; Raadsen et al., 2016, 2018; Yperman, 2007) is based on a re-formulation of the LWR equations in terms of cumulative flows. These models determine cumulative numbers of vehicle that cross link boundaries and describe the influence of traffic changes at one link boundary on the other one. They generally do not require discretization. However, they do not describe flow conditions within the link or the effect of traffic changes on the same boundary. Other approaches to solve this formulation include use of the Lax-Hopf algorithm (Mazaré et al., 2011; Simoni et al., 2020) or dynamic programing (Daganzo, 2005a; 2005b). However, these approaches are restricted to concave FDs and are not easily scalable.

Traffic states within links may be solved by tracking kinematic waves. Wong and Wong (2002) presented an analytical , for which an initial linear density profile along the link remains linear over time (Whitham, 2011). Therefore, the link exhibits piecewise linear densities. Piece endpoint locations and slopes change depending on the shockwaves generated at the link ends. Cai et al. (2009) and Lu et al. (2008, 2009) extended this approach to other concave FD shapes. Henn (2005) proposed a wave tracking solution algorithm for a network. It is based on time and space discretization and resolves interactions among waves. However, this approach is computationally expensive and only applicable with concave FDs.

This study proposes a model to solve the LWR equations based on the concept of wave tracking. However, it tracks shockwaves, the disturbances generated in the flow, rather than the kinematic waves themselves. It uses flow conservation through links to determine times of state changes at link ends with no assumption regarding the shape of the FD except that the space mean speed is a decreasing function of density. This approach may also be viewed as tracking the propagation of information on discontinuities in the flow. It can be used more generally to propagate other types of information, such as changes in routing, through the network. At nodes, flows from incoming links are distributed to outgoing ones with a procedure that captures route choices and guarantees flow conservation.

A simulation model that implements this approach is developed in a combination of time based and event-based procedure which makes the model trackable for large scale applications. The rest of this paper is organized as follows. Section 2 presents the tracking of the propagation of flow discontinuities within a link using shockwaves. Section 3 presents a node model that captures the flow change and routing procedure. Section 4 introduces the IPM concept with its basic elements. The simulation process in a distributed structure is proposed in Section 5. Section 6 examines the IPM in a various case-studies. Finally, some concluding remarks on findings and future work are presented in Section 7.

The main contribution of this paper is to introduce the information propagation concept which is a generalized approach of the wave tracking algorithms. The model tracks the information that are relevant for representing the traffic dynamics through the transportation facilities. The developed model describes the creation, propagation, interaction, and dissipation of different types of information through the network. The proposed representation allows us to cluster the data within links and generates a flexible distributed simulation model.

To the best of the authors knowledge, this is the first study that introduces the transportation network dynamics as information exchange between the networks elements, which will extend the macroscopic models' capability to represent complex systems, special attention will be given on the information exchange between toll lanes facilities.

## 2. Tracking disturbances within a link

Links are assumed to be internally homogenous. Therefore, bottlenecks and disturbances in flow are generated only at the link boundaries and then propagate through the link. The proposed model is based on determination of the time of arrival of these shockwaves to the other link boundary (or their termination within the link) and tracking of the flow states at the link boundaries.

Suppose that at a link boundary point $x$ the density and flow rate is at a steady state with $k_1(x, t)$ and $f_1 = f(k_1)$, respectively. At time $t$, a disturbance in flow occurs at the boundary, and flow conditions change to $k_2(x, t)$ and $f_2 = f(k_2)$. This change generates a shockwave that defines the boundary between the two traffic states. The shockwave propagates along either the upstream or downstream link with a speed dictated by flow conservation considerations:

$$s(t) = \frac{f_1 - f_2}{k_2 - k_1} \qquad (3)$$

Where, $s(t)$ is the shockwave speed. $k_1$ and $f_1$ are the initial density and flow rate. $k_2$ and $f_2$ are the values after the change.

Within the link, shockwaves may intersect with adjacent ones. Two or more shockwaves that intersect are terminated and may generate a new shockwave with speed that is determined by equation (3) above using the most upstream and downstream flow regimes at the intersection point. In case these two flow regimes have identical characteristics, a new shockwave will not be generated. The intersection time and location of two adjacent shockwaves are given by:

$$\tau = \frac{x_2 - x_1 + t_1 s_1 - t_2 s_2}{s_1 - s_2} \tag{4}$$

$$x = x_1 + (\tau - t_1)s_1 = x_2 + (\tau - t_2)s_2 \tag{5}$$

Where, $\tau$ is the intersection time. $x_1, t_1, s_1$ and $x_2, t_2, s_2$ are the generation position and time of the two intersecting shockwaves and their speeds, respectively.

A shockwave that does not intersect with other ones, continues to propagate until it reaches the other link boundary. There, it changes the flow state and terminates. A new shockwave may be generated because of the flow state change. Figure 1 demonstrates the shockwave propagation within a link by showing the link's time-space diagram (a) and the corresponding fundamental diagram, state transitions and shockwave speeds (b). The figure shows five shockwaves that are generated at the link boundaries at different times. The first intersection among three of them occurs at time $\tau_1$. As a result, the intersecting shockwaves are terminated. A new shockwave is not generated since both upstream and downstream streams are in state A. At time $\tau_2$, a second intersection occurs. This time, a new shockwave is created at the intersection points and propagates to the link boundary with a new speed.

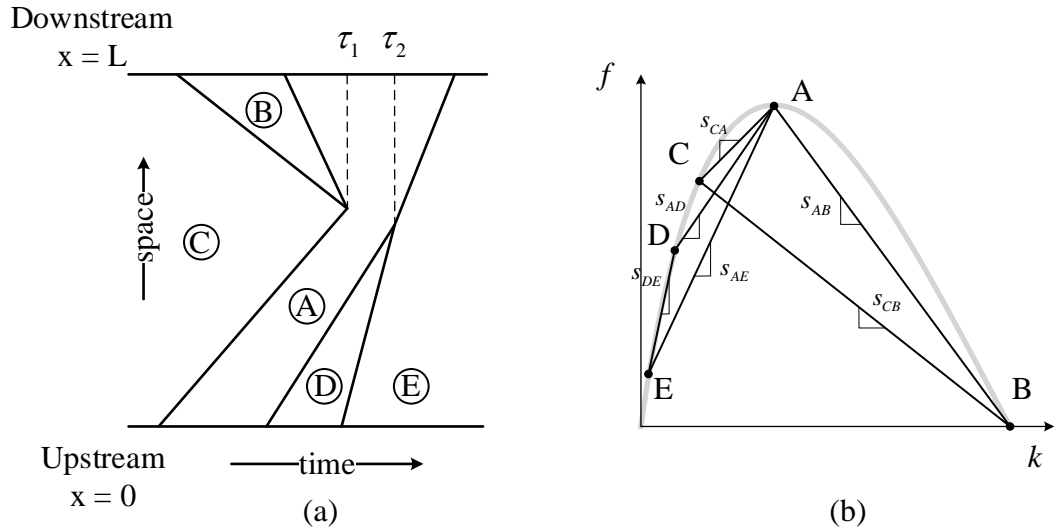

Figure 1. (a) Space-time diagram of flow distribution (b) Flow density relationship

## 3. Flow allocation at nodes

An extension of the link model above to the network level requires using nodes to represent intersections, points of geometrical changes or network boundaries. The main functionality of nodes is to allocate incoming upstream flows among the downstream links based on routing rules and considering the node and link capacities.

It is assumed that traffic demands use pre-specified sets of routes between the relevant origins and destinations. Time-dependent route proportions are calculated using a route choice model. Large scale networks may include very large numbers of routes. To simplify their representation, they are converted to a local form, in which each link maintains a list of the routes that pass through it, as shown in Figure 2. A route-route transition matrix $T_n[R \times S]$ maps the routes on the upstream links to the ones on the downstream links. Entries $T_n(r,s)$ are equal 1 if incoming route $r$ maps to outgoing route $s$, and 0 otherwise. This representation allows to merge routes that share the same links from the current node to the destination. These matrices are invariant as they depend only on the network topology and route sets. Therefore, they only need to be calculated once.

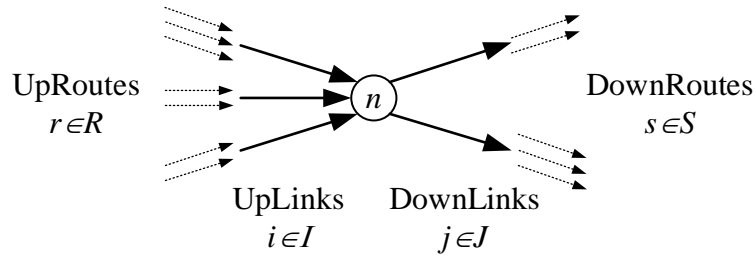

Figure 2. Representation of links and routes at a node

Given the flows on upstream links and the proportions of route flows that make them up, the flows on downstream links may be calculated, suppressing the node index, as:

$$F_J = P_{IJ} \cdot F_I = A_{JS} \cdot T \cdot diag(P_R) \cdot (A_{IR})' \cdot F_I \qquad (6)$$

Where, $F_I[I \times 1]$ and $F_J[J \times 1]$ are the exit flows on the upstream links and the entering flows on the downstream links, respectively. $P_{IJ}[I \times J]$ is a matrix of upstream link to downstream link flow proportions. $A_{JS}[J \times S]$ and $A_{IR}[I \times R]$ are link-route incidence matrices, with entries $A_{JS}(j,s)$ or $A_{IR}(i,r)$ that are equal 1 if route $s$ (or $r$) uses link $j$ (or $i$), and 0 otherwise. $P_R[R \times 1]$ is an array of the proportions of route flows within the link flows. Thus, the entries for all routes that use the same link sum up to 1. The $diag(\cdot)$ operator transforms a column array into a square matrix with the array entries on its main diagonal.

Flows through a node may be constrained by the entry capacities into the downstream links, which are determined by the flow regimes at their upstream endpoints. In under-saturated conditions, the incoming flow may reach the maximum flow allowed by the link's FD. In over-saturated conditions, the incoming flow is limited by the prevailing flow:

$$C_j = \begin{cases} f(k_j^o) & k_j \leq k_j^o \\ f(k_j) & k_j > k_j^o \end{cases} \tag{7}$$

Where, $C_j$ is the inflow capacity of link $j$. $k_j$ and $k_j^o$ are the current and optimal density on the link, respectively.

As a result, arriving flows on the upstream links may not all be able to pass through the intersection. It is assumed that the flows through the node follow the following conditions:

1. The proportions of turning movements in upstream link flow are kept. Thus, the most constrained turn from the link dictates the flows on the other movements from it. At the extreme, if a specific turn is blocked, other movements from that link would also be blocked.
2. The allocation of capacities of downstream links among the upstream links that flow into them is determined by a priority matrix $W[I \times J]$, which captures the effects of traffic signs and control. The matrix is scaled such that the entries in each column, which represent the priorities for a specific downstream link, sum up to 1. Unused turn capacities (i.e., when the upstream turn flow is lower than the allocated capacity) are re-distributed among the other upstream link using the same priorities.

The procedure presented in Pseudo-code 1 guarantees satisfying these requirements. $D_I$ are the flow demands arriving to the node from the upstream links. $D_{IJ}$ are demanded turn movements. $C_J^*$, $D_{IJ}^*$ and $F_{IJ}^*$ are the residual capacities, demands and flows, respectively. $UpList$ is a list of upstream links that may send flows to the downstream links. After initialization of these matrices (lines 2-3), downstream link capacities are allocated to the turn movements using the priority matrix $W$ (line 5). The parameters $\alpha_i$ (line 6) represent the proportion of demand that can be sent on the various turns from link $i$ without exceeding any of their capacity allocations. The use of a single value for each upstream link (line 7) guarantees that the assigned flows exiting from this link maintain the turn proportions. Next, the turn flows residual capacities and residual demands are updated, and so is the list of upstream links that may send additional flow through the node (lines 8-10). $1_I [I \times 1]$ is a unit column vector. No residual demand for a turn means that the upstream link does not have any remaining demand to send and so it is removed from $UpList$. No residual capacity on a downstream link means that it cannot receive additional flows from any of the upstream links connected to it. However, to maintain the turn proportions from these links, they cannot send flows to any other downstream link as well. Thus, they are all removed from $UpList$. The process of capacity and flow allocation is repeated until the $UpList$ is emptied, which indicates that no upstream link can send additional flows.

Pseudo-code 1. Flow allocation at a node

| | |
|---|---|
| 1 | **Input**: $D_I, P_R, C_J, W$ |
| 2 | $D_{IJ} = P_R \cdot diag(D_I)$ |
| 3 | $F_{IJ} = 0, \ C_J^* = C_J, \ D_{IJ}^* = D_{IJ}, \ UpList = \{1, \ldots, I\}$ |
| 4 | while $UpList \neq \emptyset$ |
| 5 | $\quad C_{IJ}^* = W \cdot diag(C_J^*)$ |
| 6 | $\quad \alpha_i = min\left(1, \ \frac{C_{ij}^*}{D_{ij}^*} \ j = 1, \ldots, J\right) \ \forall i = 1, \ldots, I$ |
| 7 | $\quad F_{IJ}^* = D_{IJ}^* \cdot diag(\alpha_I)$ |
| 8 | $\quad F_{IJ} = F_{IJ} + F_{IJ}^*; \ C_J^* = C_J^* - (1_I)' \cdot F_{IJ}^*; \ D_{IJ}^* = D_{IJ}^* - F_{IJ}^*;$ |
| 9 | $\quad$ if $d_{ij}^* = 0$, remove $i$ from $UpList$ |
| 10 | $\quad$ if $C_j^* = 0$, remove $\{\forall i: link \ (i,j) \ exists\}$ from $UpList$ |
| 11 | end |
| 12 | **Output**: $F_{IJ}$ |

The flows on the upstream and downstream links may be calculated the turn flows, $F_{IJ}$, as:

$$F_I = F_{IJ} \cdot 1_J \tag{8}$$

$$F_J = (1_I)' \cdot F_{IJ} \tag{9}$$

Flows and their proportions on the downstream routes are given by:

$$F_S = T \cdot diag(P_R) \cdot (A_{IR})' \cdot F_I \tag{10}$$

$$P_S = \left( diag\left( (A_{JS})' \cdot A_{JS} \cdot F_S \right) \right)^{-1} \cdot F_S \tag{11}$$

Where, $F_S[S \times 1]$ and $P_S[S \times 1]$ are the flows and proportions of the link flows on the routes within the downstream links, respectively.

## 4. Information packages and their propagation

Shockwaves may be generated not only from the traffic flow dynamics, but also from other disturbances, such as time-varying travel demand and routing, provision of travel time information, dynamic traffic control and tolling, occurrence of crashes or moving bottlenecks. The effects of these disturbances on traffic flow in the network may be modeled as a generalization of the tracking of flow shockwaves to that of information packages (IP). IPs are characterized by the information they hold. Their positions and speeds depend on the information type and on the surrounding traffic flow characteristics. As with shockwaves, IPs may intersect with other ones, leading to changes in their properties. Together with the nodes' states, the IP list characterizes the state of the network and its dynamics. Various types of information can be modeled this way, such as:

**Traffic flow characteristics**: Shockwaves may be considered an IP that carries information on the flow regime. They may be generated whenever traffic flow changes: At origin nodes when the trip demands change, at any node when another IP arrives and changes the node's flow state, or when the node capacity or routing at the node change (e.g., with new travel information). Their speeds and behavior when intersecting with other shockwaves are as described above. They are unaffected by intersection with a route choice proportions IP.

**Route choice proportions**: Route proportions within a link are captured by $P_S$. Changes in this information may occur at nodes when the flow allocation process (Pseudo-code 1) is applied. This occurs following changes in the upstream flows and route proportions, the node capacities and priorities or vehicle rerouting following reception of new travel and tolling information and guidance. The IPs that carry the routing information travel only in the downstream direction and with the space mean speed of the flow stream they travel within. This means that routing IPs cannot intersect with each other. They can intersect with traffic flow IPs (shockwaves). The flow space mean speed is always larger than the shockwave speed, and so routing IPs may intersect only with shockwaves that are downstream. When an intersection occurs, the routing IP will change its speed to the prevailing one in the new flow regime. The information it carries (i.e., route proportions) will not change. The shockwave IP will not be affected.

**Moving bottlenecks**: These IPs represent disturbances such as slow vehicles or rolling road works. They are generated externally of the traffic flow and move at the lower among their own free speed ($v_b^{free}$) and the speed of the traffic stream they travel within. A bottleneck exists at the moving position of the IP. Its capacity ($C_b$) depends on the extent that it blocks or disrupts flow across it. Figure 3 illustrates the moving bottleneck and traffic states near it on the FD. Following Munoz and Daganzo (2002), Daganzo and Laval (2005) and Simoni and Claude (2017) The moving bottleneck is considered active when the stream flowing across it is interrupted. This happens at densities that are in the range defined by the intersection points (A and B in the figure) of the FD and a straight line with an offset of the moving bottleneck capacity and a slope of the moving bottleneck speed. When the moving bottleneck is active, it functions as a shockwave with a speed that equals to the bottleneck free speed. The two points A and B define the traffic states downstream and upstream of the moving bottleneck, respectively. The information carried in the IP includes its free speed and capacity and the route it follows in the network. A moving bottleneck IP may intersect with flow shockwaves propagating either upstream or downstream. If the bottleneck is active, the traffic states downstream and upstream change to those shown in Figure 3 and propagate according to equation (3). Figure 4 shows an example of a moving bottleneck and the shockwaves it generates. The moving bottleneck is shown by the dotted line. Initially the bottleneck is within a stopped queue (state A), which also determines that the bottleneck does not move. When the head of the queue is released, a backward propagating shockwave is generated transitioning to flow state B. When this shockwave intersects with the bottleneck, it terminates, and two new shockwaves are generated: One that propagates downstream with the unconstrained flow state C and the other that propagates upstream with the constrained state B. The bottleneck IP itself starts to move at its own free speed. Shockwave IPs may also be generated when the free speed or the capacity of moving bottleneck change.

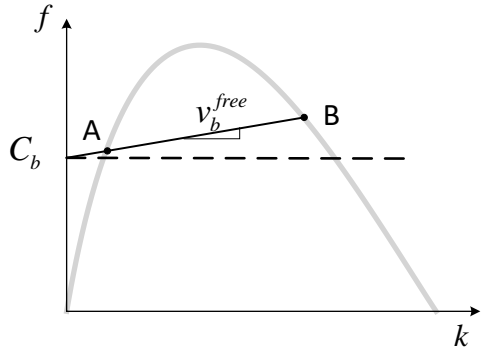

Figure 3. Traffic flow near a moving bottleneck

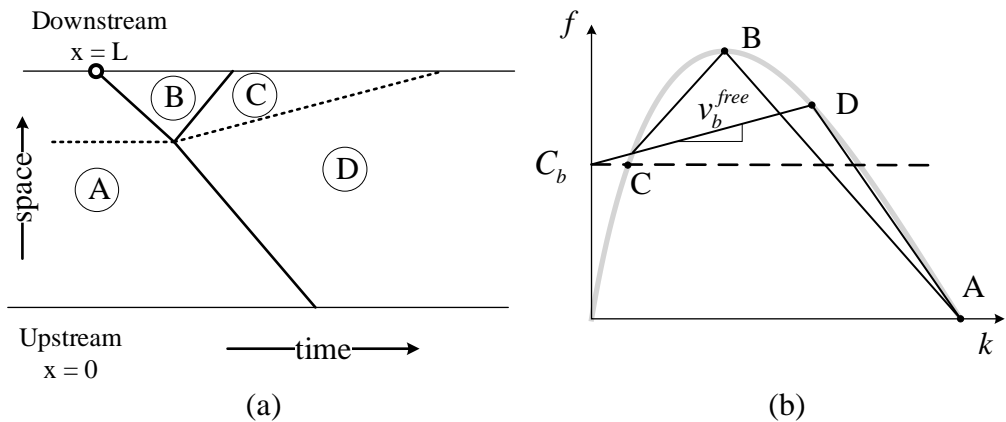

Figure 4. (a) Space-time diagram and (b) Flow density relationship near a moving bottleneck

## 5. Simulation model

### 5.1. Generic structure

A simulation model was developed that uses the IP theory outlined above to solve the dynamic state of a network. The model generates IPs, tracks them in the network and changes their properties as they interact with each other and at nodes. Figure 5 presents a flowchart of the simulation process.

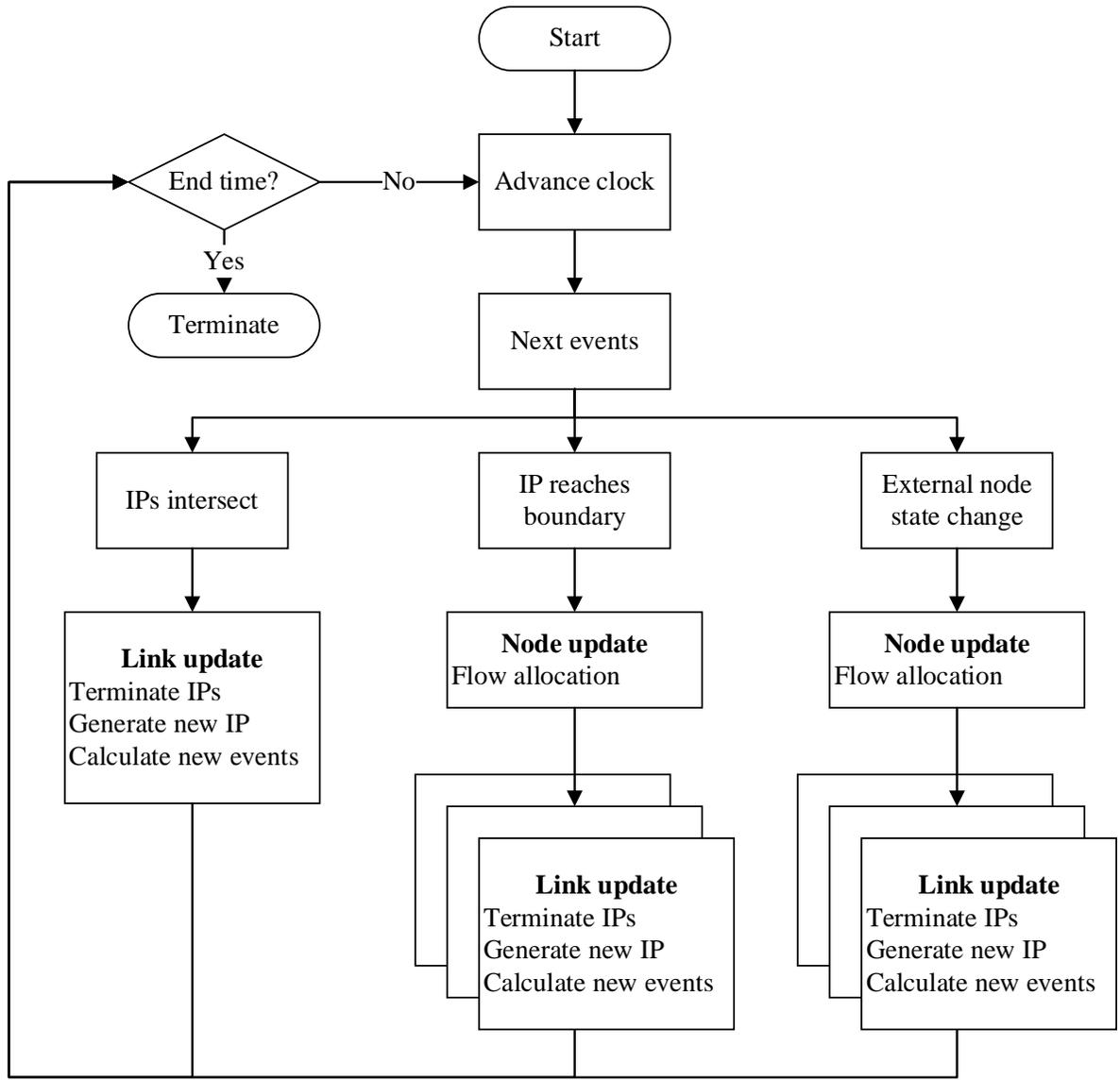

Figure 5. Flowchart of the simulation process

The model consists of two basic functions: *node update* and *link update*. The *node update* modifies nodes' states and notifies the connected links. Changes in the node state may be external to traffic flow, such as resulting from flow control and rerouting or changes in route demands at origin nodes. They may also be triggered by arrival of propagated IPs on connected links. Change in a node's state requires rebalancing of its flows and routes proportions allocation as described in Pseudo-code 1. The change in the node state may generate new IPs on connected links. Therefore, the *link update* functions are invoked for these links. The *link update* function generates the new

IPs and calculates new events related to them: IPs may intersect with each other or reach the link boundary. When IPs intersect, based on equation (4), the *link update* function terminates these IPs and generated new ones. An IP that reaches the link boundary triggers a *node update*. The figure shows a time-paced implementation, but this structure is suitable also for event-based variation. It may also be adapted for distributed representation, which is described next.

## 5.2. Distributed implementation

A distributed implementation may be useful to reduce run times for large networks. This can be achieved by parallelization of the updates of different links. Parallel processing of links requires that the IPs within it may be propagated independently of any other link, which is the case if no new IPs are generated at the link boundaries. Therefore, a time-paced two-stage processing of IPs is proposed: First, nodes process the IPs that are close to them (i.e., that may reach the link boundary within the time step), up to the point where they get further away and transferred to the link. Then, links process the IPs that are further from the nodes (i.e., cannot reach the link boundaries) and those that were transferred from the nodes.

The implementation is based on splitting links into three zones, as shown in Figure 6.: upstream and downstream node zones and the remaining middle zone. The node zones are defined such that any IP that outside of them, cannot possibly reach the link boundary within the time step. Thus, length is largest at the beginning of the time step and linearly decreases:

$$L_{n,l}(t) = v_l^{free}(t_0 + dt - t) \tag{12}$$

Where, $L_{n,l}(t)$ is the length of node $n$ zone on link $l$ at time $t_0 \leq t \leq t_0 + dt$. $v_l^{free}$ is the free (maximum) speed in the link. $t_0$ is the step start time. $dt$ is the time step size. To avoid overlap among the node zones in the same link, the step size is bounded, such that:

$$dt \leq \min_l \left(\frac{L_l}{2v_l^{free}}\right) \tag{13}$$

Where, $L_l$ is the length of the link.

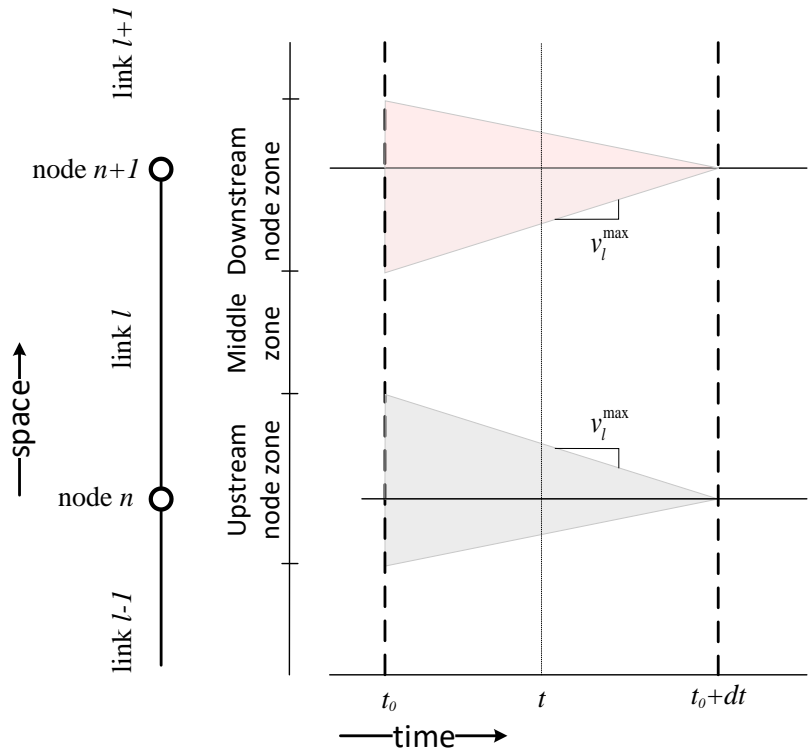

Figure 6. Time-space diagram of link zones within a time step

Pseudo-code 2 presents the overall flow of the distributed IPM. It uses two new functions: *node zone update* and *middle zone update*, which are similar to a link update but limited to their respective parts of the link. First, node zones are handled in parallel. Within these zones, IPs may reach the link boundary and trigger a *node update* or reach the boundary of the node zone. In the latter case, the IPs will not be further advanced by the node zone. The node zone exit time and location are calculated by equation

After all node zones have been processed, the middle zone updates are invoked. They consider the IPs that were within the middle zone at the start of the time step and those that exited the node zones.

Pseudo-code 2. Distributed IPM

| | |
|---|---|
| 1 | $for\ t = t_0 : dt : T$ |
| 2 | *Do in parallel* |

```
3           node zone update (n):   ∀n ∈ N
4       end
5       Do in parralel
6           middle zone update(l):   ∀l ∈ L
7       end
8   end
```

## 5.3. Extracting route travel times

The route travel time for a vehicle exiting the network can be calculated by tracking the travel times of the route's links. The time dependent travel times for each link are part of DNL outputs. Let $\boldsymbol{\tau}_l^-$ be vector of the cumulative inflows and outflows profiles for link $l$. $\Theta(\cdot)$ calculates the travel time for exiting vehicles at time $t$:

$$\tau_l^-(t) = \Theta\left(\boldsymbol{\tau}_l^-,\ t\right), \quad l \in \mathbf{L} \qquad (14)$$

Let $l_1, l_2, \ldots, l_j$ are the links' sequence for route $r$ sorted from origin to destination. presents the route travel time calculation. The algorithm starts from the last link and calculates the entry time for vehicles exited at time $t$. The entry time is the exit time for the previous link. The process is repeated recursively until it reaches the first link. It should be noted that this process can be applied during the simulation process or as a post process.

Pseudo-code 3. Calculating route travel time

| Calculating single route travel-time |
|---|
| 1 **input :** $t$, $\tau^-_l$, $l \subseteq \mathbf{L_r}$ |
| 2 $\tau^- = 0$ |
| 3 $i = j$ |
| 4 **while** $i > 0$ |
| 5      $\tau^- \leftarrow \tau^- + \tau^-_{l_i}(t)$ |
| 6      $t \leftarrow t - \tau^-$ |
| 7      $i \leftarrow i - 1$ |
| 8 **end** |
| 9 **output :** $\tau^-$ |

## 6. Case study

The IPM was applied to the Southbound direction of the Ayalon Highway (Route 20) in Tel Aviv, Israel. The Ayalon Highway, shown in Figure 7, runs on the eastern side of central Tel Aviv and connects all the major highways leading to the city. The modeled section of this highways is about 13 Km long. It includes 8 onramps and 8 offramps. Its network representation consists of 57 links, 58 nodes, 13 centroids and 19 centroid connectors. 25 sensor stations are located upstream of each on- ramp and downstream of each off-ramp. Travel demand on the highways is for 49 OD pairs. Demand matrices for the morning (6AM – 9AM) and the afternoon (3PM - 6PM) peak periods were used with 34,057 and 28,673 trips, respectively.

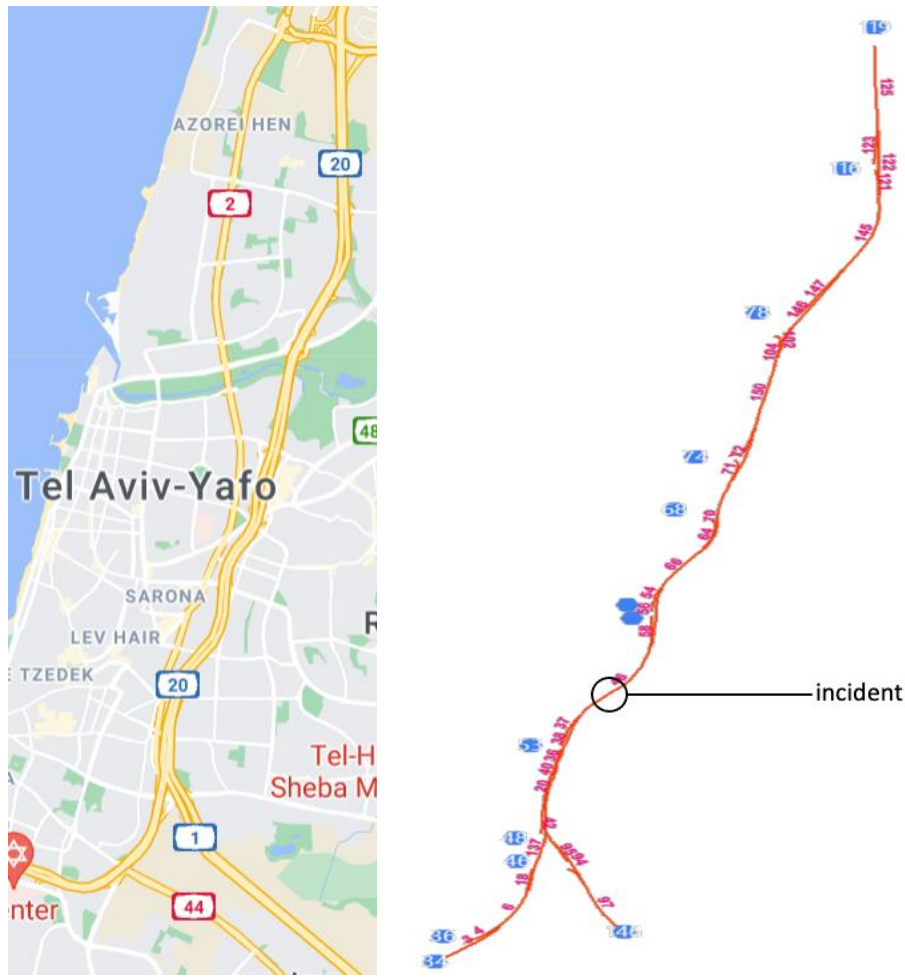

Figure 7. The Ayalon Highway network

## 6.1. Calibration

The IPM model was calibrated against a TransModeler micro-simulation model (TSM) of the same section that was previously developed and calibrated using real-world sensor data. The calibration against a micro-simulation model rather than the real-world observations is useful for testing of the model's response to various conditions for which observations may not be available in the real-world data.

The calibration involved 41 parameters including:

- Fundamental diagrams (18 parameters): Triangular FDs were assumed for each of 6 different link types. Each FD is defined by three parameters: maximum flow, critical density, and maximum density.

- Bottlenecks (15 parameters): Mainly consist of offramp boundaries and merging nodes.
- Priorities (8 parameters): Defined for each of the on-ramp merging nodes.

The calibration objective was to minimize the deviation of sensor counts and origin to destination travel times calculated for 5-minute intervals with the IPM against those measured in the TSM simulations. A weighted least squares objective was used:

$$\min_{\theta} \sum_{t=1}^{T} \left( \sum_{s=1}^{N_S} \left( y_{s,t}^{TSM} - y_{s,t}^{IPM}(\theta) \right)^2 + w_\tau \sum_{od=1}^{N_{OD}} \left( \tau_{od,t}^{TSM} - \tau_{od,t}^{IPM}(\theta) \right)^2 \right) \quad (15)$$

Where, $\theta$ are the parameters to be calibrated. $t = 1 \ldots T$ are 5-minute time intervals. $y_{s,t}^{TSM}$ and $y_{s,t}^{IPM}$ are counts at sensor station $s$ and time interval $t$ in the TSM and IPM model, respectively. $\tau_{od,t}^{TSM}, \tau_{od,t}^{IPM}$ are origin-destination travel times for OD pair $od$ departing in time interval $t$ in the TSM and IPM model, respectively. $N_S$ and $N_{OD}$ are the numbers of sensor stations and OD pairs, respectively. $w_\tau$ are relative weight assigned to the deviations in OD travel times in the objective function.

The weights $w_{OD}$ capture the inverse of the variance of the errors in OD travel time measurements, relative to those of the sensor counts. They are not known a-priori. Therefore, an iterative weighted least squares procedure (Fuller and Rao 1978, Chen and Shao 1993) was used to determine the weights and calibrate model parameters simultaneously. First, equation (15) is solved with initial weights. Then, the variances of the errors of the counts and OD travel times are estimated with the calibrated parameters and used to calculate new weights:

$$w_\tau^{(k+1)} = \frac{MSE_y^{(k)}}{MSE_\tau^{(k)}} \quad (16)$$

Where, $w_\tau^{(k+1)}$ are the weights of the OD travel time parameters in iteration $k + 1$. $MSE_y^{(k)}$ and $MSE_\tau^{(k)}$ are the mean squared deviations of the sensor counts and OD travel times, respectively, between the IPM and TSM models in iteration $k$.

The new weights are used to re-optimize equation (15). The iterative process continues until the change in weights between consecutive iterations is sufficiently small. The calibration process used data from both the AM and PM peak periods.

## 6.2. Experiments

The calibrated model was applied in three experiments using scenarios with different demands and capacity reductions (e.g., due to an incident). Demands were based on AM period which is more congested. In all cases the IPM results are compared against those of the TSM, which it was calibrated against. The experiments are summarized in

Table 1.

In experiment 1, the demand matrix is uniformly scaled by a multiplier $\theta$. This multiplier is small than a unit since the TSM network operates at its capacity with the base demand. Higher levels of demand led to gridlock situations with large queues of vehicles unable to enter the network. In experiment 2, each entry in the demand matrix was randomly perturbed such that the expected value of the total demand in the network is preserved. The perturbations are random and so, after initial testing of the stochasticity of the results, the experiment was repeated 20 times with different draws of random numbers to achieve a 95% confidence level for an interval of error of 10% of the values of the goodness of fit measures that were calculated. The perturbation is given by:

$$OD_{ir} = OD_i + 2(R_{ir} - 0.5)\alpha \qquad (17)$$

Where, $OD_{ir}$ is the demand value for OD pair $i$ in repetition $r = [1,..,10]$. $OD_i$ is the base travel demand for the OD pair. $R_{ir} \sim U(0,1)$ is a random number. $\alpha$ is a perturbation factor.

In experiment 3, the effects of incident causing capacity reductions were evaluated. The incidents were located on the mainline in a section with five lanes. The incident location is marked in Figure 7. In TSM, lanes, from right to left, were blocked according to the number of blocked lanes specified in each scenario. In IPM a comparable reduction in link exit capacity was specified at its downstream end. The experiment was run using the AM peak scenario with the incidents occurring at 8:15 AM and lasting 10 minutes.

Table 1. Preformed Experiments.

| Experiment | Factor being varied | Levels |
|---|---|---|
| 1 | OD scale | $\theta = 0.6, 0.7, 0.75, 0.8, 0.85, 0.9, 0.95, 1.0$ |
| 2 | OD structure | $\alpha = 0.05, 0.1, 0.15, 0.2$ |

| 3 | Incident (capacity reduction) | $N = 1, 2, 3, 4, 5$ lanes out of five |

## 6.3. Goodness of fit measures

Both for calibration and in the experiments, three goodness-of-fit measures were calculated: The root-mean-square error (RMSE), root-mean-square percent error (RMSPE) and Theil's inequality coefficient (U). These measures are given by:

$$RMSE = \sqrt{\frac{1}{N}\sum_{n=1}^{N}(Y_n^{TSM} - Y_n^{IPM})^2} \tag{18}$$

$$RMSPE = \sqrt{\frac{1}{N}\sum_{n=1}^{N}\left(\frac{Y_n^{TSM} - Y_n^{IPM}}{Y_n^{TSM}}\right)^2} \tag{19}$$

$$U = \frac{\sqrt{\frac{1}{N}\sum_{n=1}^{N}(Y_n^{TSM} - Y_n^{IPM})^2}}{\sqrt{\frac{1}{N}\sum_{n=1}^{N}(Y_n^{TSM})^2} + \sqrt{\frac{1}{N}\sum_{n=1}^{N}(Y_n^{IPM})^2}} \tag{20}$$

Where, $Y_n^{IPM}$ and $Y_n^{TSM}$ are the IPM simulated and TSM measurements, respectively, at space-time interval point $n$.

Theil's inequality coefficient, $U$ is bounded, $0 \leq U \leq 1$, where, $U = 0$ implies perfect fit between the observed (i.e. TSM simulated) and simulated measurements, and $U=1$ implies the worst possible fit. It can be decomposed to three proportions of inequality: $U^M$, $U^S$ and $U^C$, which capture the proportions of bias, variance, and covariance in the error. These terms sum up to a unit, and the first two proportions should be as small as possible. They are given by:

$$U^M = \frac{(\bar{Y}^{TSM} - \bar{Y}^{IPM})^2}{\frac{1}{N}\sum_{n=1}^{N}(Y_n^{TSM} - Y_n^{IPM})^2} \tag{21}$$

$$U^S = \frac{(S^{TSM} - S^{IPM})^2}{\frac{1}{N}\sum_{n=1}^{N}(Y_n^{TSM} - Y_n^{IPM})^2} \tag{22}$$

$$U^C = \frac{2(1-\rho)S^{TSM}S^{IPM}}{\frac{1}{N}\sum_{n=1}^{N}(Y_n^{TSM} - Y_n^{IPM})^2} \tag{23}$$

Where, $\bar{Y}^{TSM}$, $\bar{Y}^{IPM}$, $S^{TSM}$ and $S^{IPM}$ are the sample means and standard deviations of the TSM measurements and IPM simulations, respectively. $\rho$ is the correlation coefficient between the two sets of measurements.

## 6.4. Results

Calibration results for both AM and PM peak periods are presented in Figure 8, which plots 5-minutes calibrated IPM sensors counts and OD travel times against the corresponding TSM measurements. The figure shows a good fit with most of the points close to the 45 degrees line, which represents a perfect fit. The fit of the calibrated IPM to the TSM measurements is quantified with the goodness-of-fit statistics shown in Table 2.

Compared to the PM peak hours, the results are less accurate for the AM peak, which experiences higher congestion levels. In both time periods both the systematic bias and variance differences are very small, which is desirable.

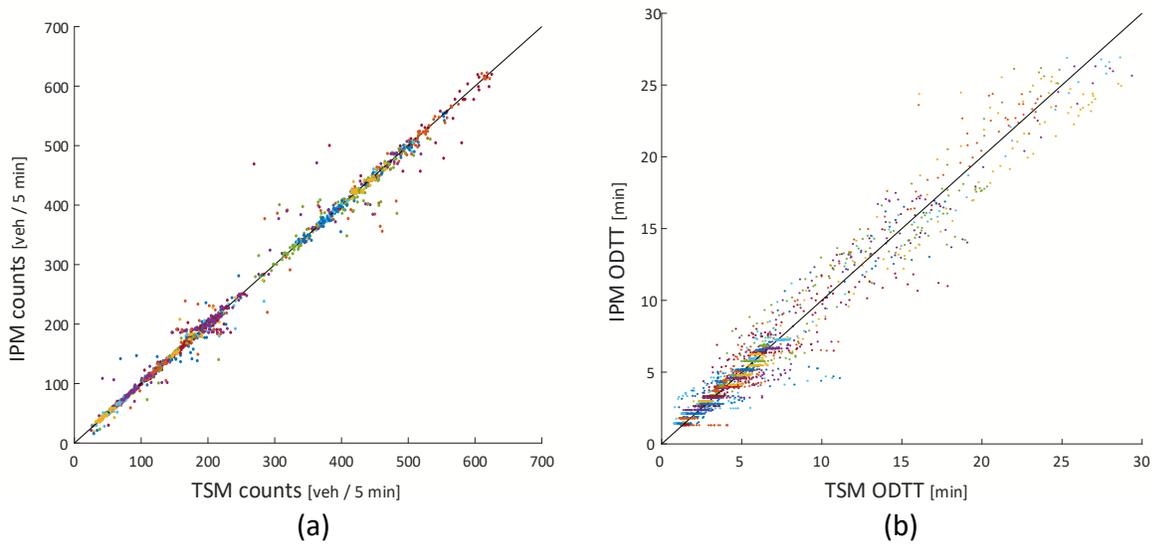

Figure 8. (a) Counts and (b) OD travel times in TSM and calibrated IPM

Table 2. Goodness-of-fit statistics for the IPM model in the AM and PM peak periods.

| Statistic | Sensor counts (veh./5 min.) | OD travel times (min.) |
| --- | --- | --- |

|   | AM peak | PM peak | AM peak | PM peak |
|---|---|---|---|---|
| $RMSE$ | 18.3 | 6.0 | 1.752 | 0.396 |
| $RMSPE$ | 0.103 | 0.036 | 0.168 | 0.064 |
| $U$ | 0.030 | 0.012 | 0.068 | 0.032 |
| $U^M$ | 0.000 | 0.003 | 0.021 | 0.007 |
| $U^S$ | 0.004 | 0.000 | 0.000 | 0.000 |
| $U^C$ | 0.997 | 0.998 | 0.980 | 0.995 |

The results of experiment 1 are presented in Table 3. The results demonstrate the inverse correlation between the congestion level and the IPM ability to replicate TSM. Higher levels of congestion lead to lower accuracy. Nevertheless, the tested scenarios with less congestion show more accurate results compared to the base scenario (calibration scenario). On the other hand, $U^m$ index for the travel times shows an opposite behavior, i.e., the systematic error increases when demand decreases. This phenomenon emerged due to shifts in the free travel times in both models. These shifts are expected since the calibration was done on congested scenarios.

Table 3. Performance measures in experiment 1

| | Demands scale $\theta$ Statistic | 1.0 | 0.95 | 0.9 | 0.85 | 0.8 | 0.75 | 0.7 | 0.6 |
|---|---|---|---|---|---|---|---|---|---|
| Counts | $RMSE$ | 18.3 | 6.5 | 5.6 | 5.4 | 5.2 | 5.0 | 4.8 | 4.4 |
| | $RMSPE$ | 0.103 | 0.045 | 0.044 | 0.046 | 0.043 | 0.044 | 0.046 | 0.044 |
| | $U$ | 0.030 | 0.011 | 0.010 | 0.010 | 0.010 | 0.010 | 0.010 | 0.010 |
| | $U^M$ | 0.000 | 0.003 | 0.013 | 0.007 | 0.006 | 0.006 | 0.010 | 0.009 |
| | $U^S$ | 0.004 | 0.001 | 0.003 | 0.000 | 0.000 | 0.000 | 0.001 | 0.000 |
| | $U^C$ | 0.997 | 0.998 | 0.985 | 0.994 | 0.995 | 0.995 | 0.990 | 0.993 |
| OD Travel times | $RMSE$ | 1.754 | 1.266 | 0.942 | 0.378 | 0.312 | 0.300 | 0.330 | 0.360 |
| | $RMSPE$ | 0.168 | 0.139 | 0.126 | 0.061 | 0.052 | 0.051 | 0.056 | 0.062 |
| | $U$ | 0.068 | 0.068 | 0.069 | 0.031 | 0.025 | 0.025 | 0.027 | 0.030 |

| | | | | | | | | |
|---|---|---|---|---|---|---|---|---|
| $U^M$ | | 0.0213 | 0.005 | 0.017 | 0.000 | 0.067 | 0.157 | 0.247 | 0.342 |
| $U^S$ | | 0.000 | 0.019 | 0.088 | 0.001 | 0.009 | 0.002 | 0.004 | 0.004 |
| $U^C$ | | 0.980 | 0.978 | 0.897 | 1.002 | 0.926 | 0.843 | 0.752 | 0.656 |

The average results over the 20 repetitions of experiment 2 are presented in Table 4. The results show a decrease in the accuracy with an increasing perturbation ~~randomness~~. This is reflected in the resulting values of the RMSE and RMSPE for both counts and travel times. This finding can be caused since the calibration was done on congested scenarios with a specific pattern. The RMSPE in OD travel times are larger than vehicle counts, since travel times are affected from the different congestion perturbations along the route.

Table 4. Performance measures in experiment 2

| | Perturbation factor $\alpha$ Statistic | 0.05 | 0.1 | 0.15 | 0.2 |
|---|---|---|---|---|---|
| Counts | $RMSE$ | 14.7 | 16.0 | 16.2 | 18.8 |
| | $RMSPE$ | 0.070 | 0.080 | 0.083 | 0.091 |
| | $U$ | 0.025 | 0.027 | 0.027 | 0.031 |
| | $U^M$ | 0.003 | 0.004 | 0.003 | 0.004 |
| | $U^S$ | 0.003 | 0.003 | 0.003 | 0.003 |
| | $U^C$ | 0.994 | 0.994 | 0.995 | 0.995 |
| OD Travel times | $RMSE$ | 1.7 | 2.0 | 2.3 | 2.4 |
| | $RMSPE$ | 0.158 | 0.178 | 0.193 | 0.215 |
| | $U$ | 0.065 | 0.078 | 0.092 | 0.094 |
| | $U^M$ | 0.057 | 0.087 | 0.109 | 0.108 |
| | $U^S$ | 0.039 | 0.052 | 0.038 | 0.029 |
| | $U^C$ | 0.906 | 0.862 | 0.855 | 0.864 |

The results of experiment 3 are presented in Table 5. The results show accuracy decrease when the number of blocked lanes increases. A small improvement in most of the measures for the counts was noted in the full block of the link. This phenomenon could be caused by the less stochasticity of the scenario.

Table 5. Performance measures in experiment 3

| | Blocked lanes $N$ Statistic | 1 | 2 | 3 | 4 | 5 |
|---|---|---|---|---|---|---|
| Counts | $RMSE$ | 16.1 | 15.5 | 19.4 | 26.6 | 24.0 |
| | $RMSPE$ | 0.066 | 0.069 | 0.091 | 0.141 | 0.139 |
| | $U$ | 0.027 | 0.028 | 0.032 | 0.044 | 0.041 |
| | $U^M$ | 0.004 | 0.004 | 0.003 | 0.005 | 0.000 |
| | $U^S$ | 0.006 | 0.005 | 0.004 | 0.005 | 0.000 |
| | $U^C$ | 0.992 | 0.992 | 0.994 | 0.991 | 1.001 |
| OD Travel times | $RMSE$ | 1.3 | 1.7 | 1.7 | 1.8 | 2.6 |
| | $RMSPE$ | 0.129 | 0.151 | 0.153 | 0.174 | 0.199 |
| | $U$ | 0.051 | 0.052 | 0.065 | 0.072 | 0.101 |
| | $U^M$ | 0.000 | 0.000 | 0.018 | 0.073 | 0.158 |
| | $U^S$ | 0.029 | 0.031 | 0.172 | 0.155 | 0.041 |
| | $U^C$ | 0.972 | 0.970 | 0.812 | 0.773 | 0.802 |

## 7. Conclusions

This paper presents a new macroscopic model for dynamic network loading. The model simulates traffic flow by tracking the propagation Ips that contain information on flow shockwaves and other disturbances through the network. It specifically calculates the times when the IPs arrive to the various nodes in the network and change their flow states. Whenever the state of a node changes, a new IP is transmitted through the relevant links to the other nodes connected to it. The IPs travel through the links. They may interact with other IPs on the link, which could change the information that they carry. When an IP reaches a node, it changes its state which leads to new transmitted IPs.

At the nodes, a node model is proposed that calculates its new state by balancing flows on the various routes passing through the node to satisfy capacity constraints, flow priorities and control actions.

Finally, it was shown that the simulation process can be computed in a parallel mode contributing to the computational efficiency. With such approach, each node can run independently within a time step and then synchronize with its neighbors.

Various scenarios were tested with the IPM in order to demonstrate its strengths and weaknesses. Results show the importance of calibration on the model prediction performance. It could be concluded that for a robust model, various network conditions should be considered while calibrating the network. The conducted tests illustrated the models' capability to simulate traffic in transportation networks under various traffic conditions. The models' accuracy indicates on their capability to serve as a tool for evaluating traffic control strategies and DTA applications.

Future work will be directed on implementing more types of data in the IPM, such as moving bottlenecks, and developing adaptive control strategies based on flowing information.